\documentclass[12pt,twoside,a4paper,openany]{article}
\usepackage[centertags]{amsmath}
\usepackage{amsfonts}
\usepackage[all]{xy}
\usepackage{graphicx}
\usepackage{amsthm}
\usepackage{newlfont}
\usepackage{amsmath}
\usepackage{eucal}
\usepackage{amssymb}
\usepackage{titlesec}
\usepackage{graphicx}
\titleformat{\subsection}[runin]% runin puts it in the same paragraph
       {\bfseries}
			 {\thesubsection}% the label and number
       {0.5em}% space between label/number and subsection title
       {}% formatting commands applied just to subsection title
       []% punctuation or other commands following subsection title

\newtheorem{thm}{Theorem}[section]

\newtheorem{prop}[thm]{Proposition}

\newtheorem{subs}[thm]{}

\begin{document}

\begin{center}
\textbf{\Large AN EFFECTIVE DIVISOR IN $\overline{M_{g}}$ DEFINED BY RAMIFICATION CONDITIONS}
\end{center}

\begin{center}
\text{Gabriel Mu\~noz}
\end{center}

\begin{abstract}
We define an effective divisor of the moduli space of stable curves $\overline{M_g}$, which is denoted $\overline{S^{2}W}$. Writing the class of $\overline{S^{2}W}$ in the Picard group of the moduli functor Pic$_{\text{fun}}(\overline{M_{g}})\otimes \mathbb{Q}$ in terms of the so-called Harer basis $\lambda,\delta_0,\ldots,\delta_{[g/2]}$, we prove that the relations among the coefficients of $\delta_1,\ldots,\delta_{[g/2]}$ are the same relations on coefficients as the Brill-Noether divisors. We present a result on effective divisors of $\overline{M_g}$ which could be useful to get the same relations on coefficients for other divisors. We also compute the coefficient of $\lambda$.
\end{abstract}

\section{Introduction}

%\addcontentsline{toc}{chapter}{Introduction}
%{\markboth{INTRODUCTION}{}} {\markright{INTRODUCTION}{}}

%{\rm\large

One of the most important recent endeavors in the field of Algebraic Geometry is to describe the birational geometry of the moduli spaces associated to curves.

Fixing the genus $g$, a topological invariant of curves, the moduli space of curves has been constructed by Mumford in the 60's; it is denoted $M_g$.
Even though Birational Geometry is concerned with general properties, 
a compactification of $M_g$ is useful. We are interested in the so-called Deligne--Mumford compactification
$\overline{M_g}$ (\cite{DM}), the compactification by adding stable curves to the boundary.

In order to try to characterize the effective cone of $\overline{M_{g}}$ and to answer other questions related to the birational geometry of $\overline{M_{g}}$, several effective divisors were computed in Pic$_{\text{fun}}(\overline{M_{g}})\otimes \mathbb{Q}$ in terms of the so-called Harer basis. The Brill--Noether divisors were computed by Harris and Mumford \cite{HMu} by the method of test curves. Also by the same method, Diaz \cite{D} and Cukierman \cite{C} computed other divisors.
Farkas computed several divisors by the same method and together with Popa \cite{FP} obtained inequalities between the first few coefficients of any effective divisor in $\overline{M_{g}}$ not contained in the boundary.

Recently, Cumino, Esteves and Gatto (\cite{CEG1},\cite{CEG2}) recomputed the Diaz and Cukierman divisors with a new approach. Instead of using test curves, the calculation was done over a general $1$-parameter family of stables curves. They used the theory of limit linear series for curves of compact type introduced by Eisenbud and Harris (\cite{EH1}), but in a slightly more general format, working for any nodal connected curves.
This approach has also been taken by Abreu \cite{A} to compute a new effective divisor in 
$\overline{M_g}$, in his thesis work under the guidance by Esteves. 

For $g=2n$, the divisor Abreu computed is defined as the closure of the locus of smooth curves $C$ having a pair of points $(P,Q)$ such that $Q$ has ramification weight at least $2$ in the linear system $H^{0}(\omega_{C}(-nP))$ and $P$ has ramification weight at least $2$ in the linear system $H^{0}(\omega_{C}(-nQ))$.
We can consider other classes of divisors which are similar to the divisor which was calculated by Abreu. 
For instance, for nonnegative integers $a,b$ such that $a+b=g$, a general problem is the calculation of the class of the divisor $R_{a,b}$ which is defined as the closure of the locus of smooth curves $C$ having a pair of points $(P,Q)$ such that $Q$ has ramification weight at least $2$ in the linear system $H^{0}(\omega_{C}(-aP))$ and $P$ has ramification weight at least $2$ in the linear system $H^{0}(\omega_{C}(-bQ))$. Notice that Abreu's thesis work addresses the case $a=b$.
A natural variant of this kind of divisors is: for each positive integer $1\leq n \leq g-2$, consider the divisor $\overline{S^{2}W_n}$ which is defined as the closure of the locus of smooth curves $C$ having a pair of points $(P,Q)$ such that $Q$ has ramification weight at least $3$ in the linear system $H^{0}(\omega_{C}(-nP))$. Our work addresses the case $n=1$. 

Thus, this work study the class in the Picard group of the functor Pic$_{\text{fun}}(\overline{M_{g}})$ of an effective divisor of $\overline{M_{g}}$. This divisor, $\overline{S^{2}W_1}$ in $\overline{M_g}$, is defined as the
closure of the locus of smooth curves $C$ having a pair of points $(P,Q)$ with $Q$ having ramification weight at least $3$ in the linear system $H^{0}(\omega_{C}(-P))$.  
Our approach is to use the method of test curves. For simplicity, we denote $\overline{S^{2}W}:=\overline{S^{2}W_1}$.

Writing the class of the divisor as 
\begin{center}
$\overline{S^{2}W}:=a \lambda-a_{0}\delta_{0}-a_{1}\delta_{1}-\ldots-a_{[g/2]}\delta_{[g/2]}$
\end{center}
and using the method of test curves,
we obtain the coefficient $a_{i}$ in terms of the coefficient $a_{1}$ for every $i>1$ and each odd integer $g\geq 5$. Our Theorem \ref{teoremaprincipal} proves the following relations:
\[a_i=(i(g-i)/(g-1))a_1, for\,\, every\,\, 2\leq i\leq [g/2].\]
Also, we compute the coefficient $a$ by using the Thom--Porteous formula and intersection theory. We get
\begin{center}
$a=9g^5-51g^4+129g^3-207g^2+174g-54$.
\end{center}

Notice that the relations in Theorem \ref{teoremaprincipal} are the same relations among the coefficients $a_1,\ldots,a_{[g/2]}$ as the Brill-Noether divisors. In order to obtain the relations, we use $[g/2]-1$ test curves, which are induced by families of flag stable curves over $\mathbb{P}^{1}$. Of crucial importance in the use of the test curves is our Proposition \ref{ecuaciones}. Our proposition is similar to \cite{HMo}, Thm $6.65$, item $2$, as both results give the same relations on coefficients in the expression for the class of effective divisors of $\overline{M_{g}}$. However, our Proposition \ref{ecuaciones} is better, as \cite{HMo}, Thm $6.65$, item $2$ is applied to divisors missing a certain locus of $\overline{M_{g}}$ parameterizing stable curves with only rational and elliptic components, but Proposition \ref{ecuaciones} is applied to divisors missing a certain subset of that locus. Moreover, that subset corresponds to certain flag curves which are very useful for degeneration arguments employing limit linear series techniques, as the rational components of the curves have at most $4$ nodes. Thus, our Proposition \ref{ecuaciones} could be used to prove that the class of a certain effective divisor has the same relations among the coefficients $a_1,\ldots,a_{[g/2]}$ as the Brill-Noether divisors. To be able to apply Proposition \ref{ecuaciones}, we will use our Proposition \ref{propflag}, which says that no flag curve can lie in $\overline{S^{2}W}$.

Our work is organized as follows: In Section \ref{prelim}, we present some preliminaries on ramification schemes, smoothings and limit linear systems, and some facts about the construction of $\overline{M_g}$ and about its associated Picard groups. In Sections \ref{curvaracional} and \ref{curvaeliptica}, we present some results on linear systems on rational and elliptic curves. Finally, in Section \ref{principal}, we compute the coefficient of $\lambda$ in the expression for $\overline{S^{2}W}$ and, by using the method of test curves, we get the relations between the coefficients of $\delta_1,\ldots,\delta_{[g/2]}$ in the expression for $\overline{S^{2}W}$.

A \textit{nodal curve} $C$ for us is a reduced, connected, projective scheme of dimension $1$ over $\mathbb{C}$ whose only singularities are nodes. The dualizing sheaf $\omega_{C}$ is an invertible sheaf over $C$ and the \textit{arithmetic genus} of $C$ is $g_{C}=h^{0}(C,\omega_{C})$.

%\newpage

%\pagenumbering{arabic} {\rm\large

%\section{Limit linear systems and ramification schemes}

%\pagenumbering{arabic} {\rm\large

\section{Preliminaries}\label{prelim}

\begin{subs}{\em(}Ramification schemes{\em)} {\em Let $\pi:\mathcal{X}\rightarrow T$ be a flat, projective morphism whose fibers are nodal curves of genus $g$. We say that $\pi$ is a \textit{family of curves}.
Suppose $\mathcal{X}$ is a nonsingular scheme.
Let $\mathcal{L}$ be an invertible sheaf on $\mathcal{X}$ and $\mathcal{V}\subseteq \pi_{*}\mathcal{L}$ a locally free subsheaf of rank $r+1$, for an integer $r\geq 0$.
Suppose for each $t\in T$ the composition
\[V_{t}:=\mathcal{V}_{t}/(\mathfrak{m}_{T,t}\mathcal{V}_{t})\longrightarrow (\pi_{*}\mathcal{L})_{t}/(\mathfrak{m}_{T,t}(\pi_{*}\mathcal{L})_{t})\longrightarrow H^{0}(\mathcal{X}_{t},\mathcal{L}\big|_{\mathcal{X}_{t}})\]
is injective. We call $\mathcal{V}$ a \textit{relative linear system}.

There exist sheaves $J_{\pi}^{i}(\mathcal{L})$ for each integer $i\geq 0$ satisfying the following properties (see \cite{E1}, \cite{LT})

(1) $J_{\pi}^{0}(\mathcal{L})\cong \mathcal{L}$.

(2) $J_{\pi}^{i}(\mathcal{L})$ is locally free of rank $i+1$.

(3) There are natural evaluation maps $e_{i}:\pi^{*}\pi_{*}\mathcal{L}\rightarrow J_{\pi}^{i}(\mathcal{L})$.

(4) For each $i\geq 1$, there is an exact sequence of truncation
\[\xymatrix{ 0\ar[r] & \omega_{\pi}^{\otimes i}\otimes \mathcal{L}\ar[r] & J_{\pi}^{i}(\mathcal{L})\ar[r]^{r_{i}} & J_{\pi}^{i-1}(\mathcal{L})\ar[r]& 0}\]
where $\omega_{\pi}$ is the relative dualizing sheaf of $\pi$. The truncation maps $r_{i}$ are compatible with the evaluation maps, i.e., $e_{i-1}=r_{i}\circ e_{i}$ for every $i\geq 1$.

When $\pi$ is a family of smooth curves, the sheaves $J_{\pi}^{i}(\mathcal{L})$ are called \textit{relative sheaves of principal parts} of order $i$ of $\mathcal{L}$.

Let $\mathcal{W}'_{\mathcal{V}}$ be the degeneracy locus of the natural evaluation map

\[u_{r}:\pi^{*}\mathcal{V}\rightarrow \pi^{*}\pi_{*}\mathcal{L}\rightarrow J_{\pi}^{r}(\mathcal{L}).\]

Notice that $u_{r}$ is a morphism between locally free sheaves of rank $r+1$ over $\mathcal{X}$. Locally, $\mathcal{W}'_{\mathcal{V}}$ is given by the zero locus of a Wronskian determinant of a sequence of $r+1$ functions. Furthermore, $\mathcal{W}'_{\mathcal{V}}$ has the property that $\mathcal{W}'_{\mathcal{V}}\cap \mathcal{X}_{t}$ is the ramification divisor $R_{V_{t}}$ of the linear system $V_{t}\subseteq H^{0}(\mathcal{X}_{t},\mathcal{L}\big|_{\mathcal{X}_{t}})$ for every smooth fiber $\mathcal{X}_{t}$.
Let $\mathcal{X}_{ns}\subseteq \mathcal{X}$ be the locus of nonsingular fibers of $\pi$. The closure $\mathcal{W}'_{\mathcal{V}}\cap \mathcal{X}_{ns}$ in $\mathcal{X}$ is denoted by
$\mathcal{W}_{\mathcal{V}}$. We call $\mathcal{W}_{\mathcal{V}}$ the \textit{ramification divisor} of $(\mathcal{V},\mathcal{L})$.

In case $\mathcal{V}=\pi_{*}\mathcal{L}$, we say that $\mathcal{W}_{\mathcal{V}}$ is the ramification divisor of the invertible sheaf $\mathcal{L}$.

Now, we will define the $k$\textit{-th special ramification locus}. The divisor $\mathcal{W}_{\mathcal{V}}$ is the zero locus of a section $w:\cal{O}_{\mathcal{X}}\rightarrow \cal{O}_{\mathcal{X}}(\mathcal{W}_{\mathcal{V}})$. By using the natural evaluation maps, this section induces derivatives
$w^{(k)}:\cal{O}_{\mathcal{X}}\rightarrow J_{\pi}^{k}(\cal{O}_{\mathcal{X}}(\mathcal{W}_{\mathcal{V}}))$. Let $S^{k}\mathcal{W}_{\mathcal{V}}$ be the zero scheme of $w^{(k)}$. We say that $S^{k}\mathcal{W}_{\mathcal{V}}$ is the $k$\textit{-th special ramification locus}. On $\mathcal{X}_{ns}$, the support of $S^{k}\mathcal{W}_{\mathcal{V}}$ is the set of points $P$ having ramification weight at least $k+1$ in the linear system $V_{\pi(P)}\subseteq H^{0}(\mathcal{L}\big|_{\mathcal{X}_{\pi(P)}})$, i.e. $wt_{V_{\pi(P)}}(P)\geq k+1$.
}
\end{subs}
\begin{subs}{\em(}Smoothings and Limit linear systems{\em)} {\em Let $C$ be a nodal curve. A \textit{smoothing} of $C$ is a flat, projective morphism $p:\cal{C}\rightarrow \Sigma$ where $\Sigma:=Spec\,\mathbb{C}[[t]]$, $\cal{C}$ is a regular scheme and $C$ is isomorphic to the special fiber. 

Let $p:\cal{C}\rightarrow \Sigma$ be a smoothing of a nodal curve $C$ of genus $g$. Let $\cal{C}_{*}$ be the generic fiber, $\mathcal{L}$ an invertible sheaf on $\cal{C}$ and $\mathcal{V}\subseteq p_{*}\mathcal{L}$ a relative linear system of rank $r+1$. Now, let $V:=H^{0}(\mathcal{V})\subseteq H^{0}(\mathcal{L})$. As $p$ is flat and $\Sigma$ is a regular, integral scheme of dimension $1$, it follows that every associated point of $\mathcal{L}$ belongs to $\cal{C}_{*}$. Then the restriction map
\begin{center}
$\Gamma(\cal{C},\mathcal{L})\rightarrow \Gamma(\cal{C}_{*},\mathcal{L}\big|_{\cal{C}_{*}})$ 
\end{center}
is injective. Indeed, suppose $s\in \Gamma(\cal{C},\mathcal{L})$ satisfies $s\big|_{\cal{C}_{*}}=0$; then we have
Supp$(s)\cap \cal{C}_{*}=\emptyset$. On the other hand, if $s\neq 0$, then we can write Supp$(s)=\overline{\{x_{1}\}}\cup\ldots \cup \overline{\{x_{m}\}}$ as a union of irreducible components and we obtain that the points $x_{1},\ldots,x_{m}$ are associated points of $\mathcal{L}$, and hence these points belong to $\cal{C}_{*}$. It follows that $s=0$.
Thus, $H^{0}(\mathcal{L})$ is a torsion-free $\mathbb{C}[[t]]$-module and hence free. Also, it follows that $V$ is a free $\mathbb{C}[[t]]$-module. Notice that, since cohomology commutes with flat base change, we have the isomorphism $H^{0}(\mathcal{L})\otimes_{\mathbb{C}[[t]]}\mathbb{C}((t))\cong H^{0}(\mathcal{L}\big|_{\cal{C}_{*}})$.
Now, let $V_{*}:=V\otimes_{\mathbb{C}[[t]]}\mathbb{C}((t))$. Since $\mathcal{V}\subseteq p_{*}\mathcal{L}$ is a relative linear system, we have the injective map
$V/tV\hookrightarrow H^{0}(\mathcal{L})/tH^{0}(\mathcal{L})$, and hence $V=V_{*}\cap H^{0}(\mathcal{L})$.

Now, given a divisor $D$ on $\cal{C}$ with support in $C$, let $V(D)_{*}$ be the image of $V_{*}$ under the natural isomorphism
$H^{0}(\mathcal{L}\big|_{\cal{C}_{*}})\cong H^{0}(\mathcal{L}(D)\big|_{\cal{C}_{*}})$. We define $V(D):=V(D)_{*}\cap H^{0}(\mathcal{L}(D))$.

If $D$ is an effective divisor on $\cal{C}$, we define $V(-D):=V\cap H^{0}(\mathcal{L}(-D))$. Also, if $D\subseteq C$ is a subcurve, define $V\big|_{D}$ as the image of $V$under the restriction map $H^{0}(\mathcal{L})\rightarrow H^{0}(\mathcal{L}\big|_{D})$. Notice that, if $D$ is an effective divisor of $\cal{C}$ and $E$ is a subcurve of $C$ such that $D$ and $E$ have no common components, then $V(-D)\big|_{E}\subseteq V\big|_{E}(-D\cdot E)$.

Let $C_{1},\ldots,C_{n}$ be the irreducible components of $C$. Since $C$ is connected, for each $i=1,\ldots,n$ there exists an invertible sheaf $\mathcal{L}_{i}$ on $\cal{C}$ of the form
\[ \mathcal{L}_{i}=\mathcal{L}(\sum_{i=1}^{n}a_{i,l}C_{l})=\mathcal{L}\otimes \cal{O}_{\cal{C}}(\sum_{i=1}^{n}a_{i,l}C_{l})\]
such that the restriction map
\begin{center}
$H^{0}(C,\mathcal{L}_{i}\big|_{C})\rightarrow H^{0}(C_{i},\mathcal{L}_{i}\big|_{C_{i}})$
\end{center}
is injective. We say that $\mathcal{L}_{i}$ has focus on $C_{i}$.
Let $V_{i}:=V(\sum_{i=1}^{n}a_{i,l}C_{l})$ and let $\bar{V_{i}}$ be the image of $V_{i}$ under the restriction map
\begin{center}
$H^{0}(\cal{C},\mathcal{L}_{i})\rightarrow H^{0}(C_{i},\mathcal{L}_{i}\big|_{C_{i}})$.
\end{center}
The dimension of $\bar{V_{i}}$ is $r+1$. We say that $(\bar{V_{i}},\mathcal{L}_{i}\big|_{C_{i}})$ is a \textit{limit linear system} on $C_{i}$.

Let $R_{i}$ be the ramification divisor of $(\bar{V_{i}},\mathcal{L}_{i}\big|_{C_{i}})$ and $W$ the ramification divisor of $\mathcal{V}$. Then (see \cite{E2})
\begin{equation}\label{limitramif}
W\cap C=\sum_{i=1}^{n}R_{i}+\sum_{i<j}\sum_{P\in C_{i}\cap C_{j}}(r+1)(r-l_{i,j})P.
\end{equation}
where $l_{i,j}:=a_{i,j}-a_{i,i}+a_{j,i}-a_{j,j}$. We call $W\cap C$ the \textit{limit ramification divisor} of $(\mathcal{V},\mathcal{L})$, and $l_{i,j}$ is called the \textit{connecting number} between $\mathcal{L}_i$ and $\mathcal{L}_j$ with respect to $C_i$ and $C_j$.
}
\end{subs}
%\begin{center}
%$q_1f_1+\ldots+q_{r+1}tg_m=0$ for some $q_1,\ldots,q_{r+1}\in \mathbb{C}[[t]]$, not all zero.
%\end{center}
%We may assume that some $q_j$ is not in the ideal $(t)$. Since the images of $f_1,\ldots,f_{r+1-m},tg_1,\ldots,tg_m$ in $V(-E)/tV(-E)$ are $\mathbb{C}$-linearly 
%independent, taking the quotient of the equality above modulo $tV(-E)$, we get that $q_j\in (t)$ for every $j=1,\ldots,r+1$, a contradiction. 

%\chapter{The moduli space of stable curves and intersection theory}
\begin{subs}{\em(}The moduli space of stable curves $\overline{M_{g}}${\em)} {\em Let $g\geq 2$ be an integer. Let $\overline{M_{g}}$ denote the coarse moduli space of stable curves. We will recall how $\overline{M_{g}}$ is constructed. Given a
Deligne-Mumford stable curve $X$, we have that $\omega_{X}^{\otimes n}$ is very ample for each $n\geq 3$. Then, we may view $X$ as a closed subscheme of degree $2n(g-1)$ of
$\mathbb{P}^{N}$, where $N=(2n-1)(g-1)-1$, as by Riemann-Roch, we have that
$h^{0}(X,\omega_{X}^{\otimes n})=(2n-1)(g-1)$ for each $n\geq 2$.

We have that $\omega_{X}^{\otimes n}\cong \cal{O}_{X}(1)$; we call such a stable curve $n$-canonically embedded. Let $H$ be the Hilbert scheme parametrizing subschemes of $\mathbb{P}^{N}$ with
Hilbert polynomial $2n(g-1)T+1-g$, and $\cal{U}\subseteq \mathbb{P}^{N}\times H$ the universal closed
subscheme. There is a locally closed subscheme $K\subseteq H$ parametrizing
$n$-canonically embedded stable curves of genus $g$. Let $\cal{V}:=\cal{U}_{K}\subseteq \mathbb{P}^{N}\times K$ be the induced subscheme and $v:\cal{V}\rightarrow K$ the family induced by the second projection $\mathbb{P}^{N}\times K\rightarrow K$. We have that $K$ is smooth (see \cite{HMo}, lemma $3.35$) and the family $v:\cal{V}\rightarrow K$ is versal. 

The group of automorphisms $PGL(N)$ of $\mathbb{P}^{N}$ acts naturally on $H$. Then, there is an induced action $PGL(N)\times K\rightarrow K$. Gieseker \cite{G} constructs $\overline{M_{g}}$ as a geometric GIT quotient of $K$ under
this action for any $n$ sufficiently large. We have that the quotient map $\Phi: K\rightarrow \overline{M_{g}}$ is
also the map induced by the family $v:\cal{V}\rightarrow K$.

Now, we will recall some facts about the Picard group of $\overline{M_{g}}$. Let $A^{1}(\overline{M_{g}})$ be its Chow group of
codimension$-1$ cycle classes and Pic$(\overline{M_{g}})$ its Picard group.  We have an isomorphism
\[A^{1}(\overline{M_{g}})\otimes \mathbb{Q}\rightarrow \text{Pic}(\overline{M_{g}})\otimes \mathbb{Q}.\]
On the other hand, we have another Picard group associated to $\overline{M_{g}}$, which is called {\em Picard group of the moduli functor} $\text{Pic}_{\text{fun}}(\overline{M_{g}})$. Roughly speaking, an element $\gamma\in \text{Pic}_{\text{fun}}(\overline{M_{g}})$ is a collection of classes $\gamma_{\pi}\in \text{Pic}(S)$ for each family of stable curves $\pi:\cal{C}\rightarrow S$, such that for each Cartesian diagram
\[\xymatrix{\cal{C}'\ar[r]\ar[d]_{\pi'} &\cal{C}\ar[d]^{\pi}\\
S'\ar[r]^{f}&S}\] 
we have $\gamma_{\pi'}\cong f^{*}(\gamma_{\pi})$.

We have an isomorphism (see \cite{HMo}, Proposition $3.88$)
\[\text{Pic}(\overline{M_{g}})\otimes \mathbb{Q}\rightarrow \text{Pic}_{\text{fun}}(\overline{M_{g}})\otimes \mathbb{Q}.\]

Also, there is an isomorphism (see \cite{HMu}, p.$50$)
\[\text{Pic}_{\text{fun}}(\overline{M_{g}})\rightarrow \text{Pic}(K)^{PGL(N)},\]
where Pic$(K)^{PGL(N)}\subseteq$ Pic$(K)$ is the invariant subgroup under the action of $PGL(N)$.

We have tautological and boundary classes in $\text{Pic}_{\text{fun}}(\overline{M_{g}})$.
There is a natural element $\lambda\in$ Pic$_{\text{fun}}(\overline{M_{g}})$, which is called a \textit{tautological class}. Given a family $\pi:\cal{C}\rightarrow S$ of stable curves, define $\lambda_{\pi}:=$det$(\pi_{*}(\omega_{\pi}))$, where $\omega_{\pi}$ is the dualizing sheaf of $\pi$.

On the other hand, for each $i=0,\ldots,[g/2]$, we define the subsets $\Delta'_{i}\subseteq K$ as follows:
$\Delta'_0$ is the set of
points $s\in K$ such that the fiber $\cal{V}_{s}$ has a connecting node, and $\Delta'_{i}$, for $i\geq 1$
is the set of points $s\in K$ such that the fiber $\cal{V}_{s}$ has a disconnecting node $P$, and the
closure in $\cal{V}_{s}$ of one of the connected components of $\cal{V}_{s}-\{P\}$ has arithmetic
genus $i$. The subsets $\Delta'_{i}\subseteq K$ are closed subsets of $K$ of codimension 1. We give them their
reduced induced scheme structures. Thus, they are Cartier divisors, because $K$ is smooth. The invertible sheaves associated to the $\Delta'_{i}$ are invariant under the action of $PGL(N)$. Let $\delta_0,\ldots,\delta_{[g/2]}$ denote the corresponding elements of Pic$_{\text{fun}}(\overline{M_{g}})$. These elements are called \textit{boundary classes}. We can also view $\lambda$ and the $\delta_i$ as elements of Pic$(\overline{M_g})\otimes \mathbb{Q}$. 

The group Pic$_{\text{fun}}(\overline{M_{g}})$ is freely generated by $\lambda$ and the $\delta_i$ for $g\geq 3$ (see \cite{AC}).
For calculations, it is useful the fact that a class $\gamma\in \text{Pic}_{\text{fun}}(\overline{M_{g}})\otimes \mathbb{Q}$ is defined by its value $\gamma_{\pi}\in \text{Pic}(S)\otimes \mathbb{Q}$ on $1$-parameter families $\pi:\cal{C}\rightarrow S$, where $\cal{C}$ is smooth. Moreover, it is enough to consider just a sufficiently general family. 
}
\end{subs}

%\chapter{Linear systems on rational and elliptic curves}

\section{Preliminaries on rational curves}\label{curvaracional}

\begin{prop}\label{rational}
Let $R_1,\ldots,R_n$ be distinct points on $\mathbb{P}^{1}$, and $a_1,\ldots,a_n$ positive integers.
Define the linear system 
\begin{center}
$V:=H^{0}(\omega_{\mathbb{P}^{1}}((a_{1}+1)R_1))+ \ldots + H^{0}(\omega_{\mathbb{P}^{1}}((a_{n}+1)R_n))$

$\subseteq H^{0}(\omega_{\mathbb{P}^{1}}((a_{1}+1)R_{1}+\ldots +(a_n+1)R_n))$.
\end{center}
Then $V$ is $(a_1+\ldots+a_n)$-dimensional and has no ramification points on $\mathbb{P}^{1}-\{R_1,\ldots,R_n\}$. Furthermore, for each $i$, the orders of vanishing at $R_i$ of the sections in $V$ are
\begin{center}
$0,\ldots,a_i-1,a_i+1,\ldots,a_1+\ldots+a_n$,
\end{center}
and the ramification weight of $V$ at $R_i$ is $wt_{V}(R_i)=\sum\limits_{j\neq i} a_{j}$.
\end{prop}
{\em Proof.} Let $\mathcal{L}:=\omega_{\mathbb{P}^{1}}((a_1+1)R_1+\ldots +(a_n+1)R_n)$. Since for each $i$
\[U_i:=H^{0}(\omega_{\mathbb{P}^{1}}((a_i+1)R_i))\cap \sum_{j\neq i}H^{0}(\omega_{\mathbb{P}^{1}}((a_j+1)R_j))\]
is contained in $H^{0}(\mathcal{L}(-\sum\limits_{j\neq i} (a_j+1)R_j))$ and $H^{0}(\mathcal{L}(-(a_i+1)R_i))$, we get 
\[U_i\subseteq H^{0}(\mathcal{L}(-(a_1+1)R_1-\ldots-(a_n+1)R_n))=H^{0}(\omega_{\mathbb{P}^{1}})=0\]
for every $i$, so the dimension of $V$ is $a_1+\ldots+a_n$. On the other hand, all complete linear systems on $\mathbb{P}^{1}$ have no ramification points, so the statement of the proposition is true if $n=1$. Suppose $n\geq 2$ and let us argue by induction on $n$.
For every $0\leq m\leq a_1+1$
\begin{equation*}
\begin{split}
V(-mR_1)&=H^{0}(\omega_{\mathbb{P}^{1}}((a_1+1-m)R_1))\oplus H^{0}(\omega_{\mathbb{P}^{1}}((a_2+1)R_2))\oplus \ldots \\
&\ldots \oplus H^{0}(\omega_{\mathbb{P}^{1}}((a_n+1)R_n)).\\
\end{split}
\end{equation*}
It follows that 
\begin{align}
\text{dim}_{\mathbb{C}}V(-mR_1)&=a_1-m+a_2+\ldots+a_n \,\,\text{for every}\,\, 0\leq m\leq a_1 \label{rational1}\\
\text{and}\,V(-a_1R_1)&=V(-(a_1+1)R_1). \label{rational2}
\end{align}
Now, consider the linear system
\begin{center}
$V':=H^{0}(\omega_{\mathbb{P}^{1}}((a_2+1)R_2))\oplus \ldots \oplus H^{0}(\omega_{\mathbb{P}^{1}}((a_n+1)R_n))$

$\subseteq H^{0}(\omega_{\mathbb{P}^{1}}((a_2+1)R_2+\ldots +(a_n+1)R_n))$.
\end{center}
Since by induction $V'$ has no ramifications points on $\mathbb{P}^{1}-\{R_2,\ldots,R_n\}$, and since $V'$ has dimension $a_2+\ldots+a_n$ and $V'(-\alpha R_1)=V(-(a_1+1+\alpha)R_1)$ for every integer $\alpha\geq 0$, it follows that 
\begin{equation}\label{rational3}
V(-(a_1+1+a_2+\ldots+a_n)R_1)=0. 
\end{equation}
Then, it follows from (\ref{rational1}), (\ref{rational2}) and (\ref{rational3}) that, the orders of vanishing at $R_1$ of the sections in $V$ are 
\begin{center}
$0,\ldots,a_1-1,a_1+1,\ldots,a_1+\ldots+a_n,$ 
\end{center}
whence $wt_{V}(R_1)=a_2+\ldots+a_n$. Analogously, for each $i$, the orders of vanishing at $R_i$ of the sections in $V$ are 
\begin{center}
$0,\ldots,a_i-1,a_i+1,\ldots,a_1+\ldots+a_n,$ 
\end{center}
whence $wt_{V}(R_i)=\sum\limits_{j\neq i} a_{j}$. Then
$wt_{V}(R_1)+\ldots+wt_{V}(R_n)=(n-1)(a_1+\ldots+a_n).$
On the other hand, since 
\begin{center}
deg$(\mathcal{L})=a_1+\ldots+a_n+n-2$ and dim$_{\mathbb{C}}V=a_1+\ldots+a_n$,
\end{center}
we have that, by Plücker formula,
deg$(R_{V})=(n-1)(a_1+\ldots+a_n)$.
Therefore, we have no other ramification points.   \hfill $\Box$

\begin{prop}\label{rationalr}
Let $R_1,\ldots,R_n$ be distinct points on $\mathbb{P}^{1}$, and $a_1,\ldots,a_n$ positive integers. 

Let $\mathcal{L}:=\omega_{\mathbb{P}^{1}}((a_{1}+1)R_{1}+\ldots +(a_n+1)R_n))$.
Define the linear system 
\begin{center}
$V:=H^{0}(\omega_{\mathbb{P}^{1}}((a_{1}+1)R_1))\oplus \ldots \oplus H^{0}(\omega_{\mathbb{P}^{1}}((a_{n}+1)R_n))\subseteq H^{0}(\mathcal{L})$.
\end{center}
Let $V_1\subseteq H^{0}(\mathcal{L})$ be a linear system of dimension $a_1+\ldots+a_n-1$ contained in $V$ and containing 
\begin{center}
$H^{0}(\omega_{\mathbb{P}^{1}}((a_{1}-1)R_1))\oplus H^{0}(\omega_{\mathbb{P}^{1}}((a_{2}+1)R_2))\oplus \ldots
\oplus H^{0}(\omega_{\mathbb{P}^{1}}((a_{n}+1)R_n))$.
\end{center}
Then either $V_1$ has no ramification points on $\mathbb{P}^{1}-\{R_1,\ldots,R_n\}$ or $V_1$ has exactly one ramification point there and the ramification is simple. Furthermore, 
\begin{center}
$wt_{V_1}(R_1)=(\sum\limits_{j\neq 1} a_{j})+a_1+\ldots+a_n-2+\epsilon_1$, where $\epsilon_1\in \{0,1\}$,
\end{center}
and for each $i\neq 1$
\begin{center}
$wt_{V_1}(R_i)=(\sum\limits_{j\neq i} a_{j})-1+\epsilon_i$, where $\epsilon_i\in \{0,1\}$.
\end{center}
\end{prop}
{\em Proof.} If $a_1=1$, then by dimension considerations
\begin{center}
$V_1=H^{0}(\omega_{\mathbb{P}^{1}}((a_{2}+1)R_2))\oplus \ldots \oplus H^{0}(\omega_{\mathbb{P}^{1}}((a_{n}+1)R_n))$.
\end{center}
It follows from Proposition \ref{rational} that, for each $i\neq 1$,
\begin{center}
$wt_{V_1}(R_i)=\sum\limits_{j\neq 1,i} a_{j}=(\sum\limits_{j\neq i}a_j)-1+\epsilon_i$, with $\epsilon_i=0$
\end{center}
and 
\begin{equation*}
wt_{V_1}(R_1)=2(a_2+\ldots+a_n)=(\sum_{j\neq 1} a_{j})+a_1+\ldots+a_n-2+\epsilon_1, \text{with}\,\, \epsilon_1=1.
\end{equation*} 

Now, assume $a_1\geq 2$. By Proposition \ref{rational}, the orders of vanishing at $R_1$ of the sections in $V$ are 
\begin{center}
$0,\ldots,a_1-1,a_1+1,\ldots,a_1+\ldots+a_n$.
\end{center}
It follows that the orders of vanishing at $R_1$ of the sections in $V_1$ are of the form 
\begin{center}
$\{0,\ldots,a_1-1,a_1+1,\ldots,a_1+\ldots+a_n\}-\{l\}$, for some integer $l$.
\end{center}
Notice that, by hypothesis,
\begin{center}
$V(-2R_1)=H^{0}(\omega_{\mathbb{P}^{1}}((a_{1}-1)R_1))\oplus H^{0}(\omega_{\mathbb{P}^{1}}((a_2+1)R_2))\oplus \ldots$
$\ldots \oplus H^{0}(\omega_{\mathbb{P}^{1}}((a_{n}+1)R_n))\subseteq V_1$,
\end{center}
and since $V_1\subseteq V$, we have $V_1(-2R_1)=V_1\cap V(-2R_1)=V(-2R_1)$. Then, by Proposition \ref{rational}, dim$_{\mathbb{C}}V_1(-2R_1)=a_1+\ldots+a_n-2$ and hence $l\leq 1$. Therefore
\begin{equation*}
wt_{V_1}(R_1)=wt_{V}(R_1)+a_1+\ldots+a_n-1-l=(\sum_{j\neq 1} a_{j})+a_1+\ldots+a_n-2+\epsilon_1,
\end{equation*}
with $\epsilon_1=1-l\in \{0,1\}$.

To show the equalities $wt_{V_1}(R_i)=(\sum\limits_{j\neq i} a_{j})-1+\epsilon_i$, where $\epsilon_i\in \{0,1\}$ and $i\neq 1$, it is enough to consider the case $i=2$. Notice that
\begin{equation*}
\begin{split}
V(-(a_2+1)R_2)&=H^{0}(\omega_{\mathbb{P}^{1}}((a_1+1)R_1))\oplus H^{0}(\omega_{\mathbb{P}^{1}}((a_3+1)R_3))\oplus \ldots \\
&\ldots \oplus H^{0}(\omega_{\mathbb{P}^{1}}((a_n+1)R_n)).\\
\end{split}
\end{equation*}
Now, consider the linear system 
\begin{equation*}
\begin{split}
V'&:=H^{0}(\omega_{\mathbb{P}^{1}}((a_1+1)R_1))\oplus H^{0}(\omega_{\mathbb{P}^{1}}((a_3+1)R_3))\oplus \ldots \\
&\ldots \oplus H^{0}(\omega_{\mathbb{P}^{1}}((a_n+1)R_n))\\
&\subseteq H^{0}(\omega_{\mathbb{P}^{1}}((a_{1}+1)R_{1}+(a_3+1)R_3\ldots +(a_n+1)R_n)).
\end{split}
\end{equation*}
It follows from Proposition \ref{rational} that
\begin{equation*}
\begin{split}
V'(-2R_1)&=H^{0}(\omega_{\mathbb{P}^{1}}((a_1-1)R_1))\oplus H^{0}(\omega_{\mathbb{P}^{1}}((a_3+1)R_3))\oplus \ldots \\
&\ldots \oplus H^{0}(\omega_{\mathbb{P}^{1}}((a_n+1)R_n))\\
&\subseteq H^{0}(\omega_{\mathbb{P}^{1}}((a_{1}-1)R_{1}+(a_3+1)R_3\ldots +(a_n+1)R_n))
\end{split}
\end{equation*}
has no ramification points on $\mathbb{P}^{1}-\{R_1,R_3,\ldots,R_n\}$ and $V'(-2R_1)$ has dimension equal to $a_1+a_3+\ldots+a_n-2$. Then $V'(-2R_1-(a_1+a_3+\ldots+a_n-2)R_2)=0$ and hence $V(-2R_1-(a_1+\ldots+a_n-1)R_2)=0$. By Proposition \ref{rational}, dim$_{\mathbb{C}}V(-(a_1+\ldots+a_n-1)R_2)=2$; then, by dimension considerations
\begin{center}
$V=V(-2R_1)\oplus V(-(a_1+\ldots+a_n-1)R_2)$.
\end{center}
Therefore $V(-(a_1+\ldots+a_n-1)R_2)\nsubseteq V_{1}$ and we get 
dim$_{\mathbb{C}}V_{1}(-(a_1+\ldots+a_n-1)R_2)=1$.
On the other hand, by Proposition \ref{rational}, the orders of vanishing at $R_2$ of the sections in $V$ are
\begin{center}
$0,\ldots,a_2-1,a_2+1,\ldots,a_1+\ldots+a_n$.
\end{center}
So the orders of vanishing at $R_2$ of the sections in $V_{1}$ are of the form
\begin{center}
$\{0,\ldots,a_2-1,a_2+1,\ldots,a_1+\ldots+a_n\}-\{l\}$,
\end{center}
for some integer $l$. Since we have dim$_{\mathbb{C}}V_{1}(-(a_1+\ldots+a_n-1)R_2)=1$, it follows that $l\geq a_1+\ldots+a_n-1$. Thus 
\begin{equation*}
wt_{V_{1}}(R_2)=wt_{V}(R_2)+a_1+\ldots+a_n-1-l=a_1+a_3+\ldots+a_n-1+\epsilon_2,
\end{equation*}
where $\epsilon_2:=a_1+\ldots+a_n-l\in \{0,1\}$. 

Finally, we will prove the first statement of the proposition. Using the equalities we have shown, we get
\begin{center}
$\sum wt_{V_1}(R_i)=(a_1+\ldots+a_n-1)n-1+\sum \epsilon_i$.
\end{center}
On the other hand, by Plücker formula, deg$(R_{V_1})=(a_1+\ldots+a_n-1)n$. Then $0\leq \sum \epsilon_i \leq 1$ and $V_1$ has $1-\sum \epsilon_i$ ramification points on $\mathbb{P}^{1}-\{R_1,\ldots,R_n\}$, counted with their respective weights. This proves the first statement of the proposition.   \hfill $\Box$

\begin{prop}\label{rationalp}
Let $R_1,\ldots,R_n$ be distinct points on $\mathbb{P}^{1}$, and $a_1,\ldots,a_n$ positive integers. 
Define the linear system 
\begin{center}
$V:=H^{0}(\omega_{\mathbb{P}^{1}}((a_{1}+1)R_1))\oplus \ldots \oplus H^{0}(\omega_{\mathbb{P}^{1}}((a_{n}+1)R_n))$

$\subseteq H^{0}(\omega_{\mathbb{P}^{1}}((a_{1}+1)R_{1}+\ldots +(a_n+1)R_n))$.
\end{center}
Let $P\in \mathbb{P}^{1}-\{R_1,\ldots,R_n\}$, and consider the linear system 
\begin{center}
$V_1:=V(-P)\subseteq H^{0}(\omega_{\mathbb{P}^{1}}(-P+(a_{1}+1)R_{1}+\ldots +(a_n+1)R_n))$.
\end{center}
Then either $V_1$ has no ramification points on $\mathbb{P}^{1}-\{R_1,\ldots,R_n\}$ or $V_1$ has exactly one ramification point there and the ramification is simple. Furthermore,
for each $i$
\begin{center}
$wt_{V_1}(R_i)=(\sum\limits_{j\neq i} a_{j})-1+\epsilon_i$, where $\epsilon_i\in \{0,1\}$.
\end{center}
\end{prop}
{\em Proof.} All complete linear systems on $\mathbb{P}^{1}$ have no ramification points, so the first statement of the proposition is true if $n=1$. Suppose $n\geq 2$ and let us argue by induction on $n$. By Proposition \ref{rational}, $V_1$ has dimension $a_1+\ldots+a_n-1$. On the other hand, we have
\begin{equation*}
\begin{split}
V(-a_1R_1)&=V(-(a_1+1)R_1)\\
&=H^{0}(\omega_{\mathbb{P}^{1}}((a_{2}+1)R_2))\oplus \ldots \oplus H^{0}(\omega_{\mathbb{P}^{1}}((a_{n}+1)R_n)),
\end{split}
\end{equation*}
then 
\begin{equation*}
\begin{split}
V_1(-a_1R_1)&=V_1(-(a_1+1)R_1)\\
&=(H^{0}(\omega_{\mathbb{P}^{1}}((a_{2}+1)R_2))\oplus \ldots \oplus H^{0}(\omega_{\mathbb{P}^{1}}((a_{n}+1)R_n)))(-P)\\
&\subseteq H^{0}(\omega_{\mathbb{P}^{1}}(-P+(a_2+1)R_2+\ldots+(a_n+1)R_n)).
\end{split}
\end{equation*}
Now, consider the linear system
\begin{equation*}
\begin{split}
V'&:=(H^{0}(\omega_{\mathbb{P}^{1}}((a_{2}+1)R_2))\oplus \ldots \oplus H^{0}(\omega_{\mathbb{P}^{1}}((a_{n}+1)R_n)))(-P)\\
&\subseteq H^{0}(\omega_{\mathbb{P}^{1}}(-P+(a_2+1)R_2+\ldots+(a_n+1)R_n)).
\end{split}
\end{equation*}
Notice that $V'=V_1(-a_1R_1)$. By Proposition \ref{rational}, we have dim$_{\mathbb{C}}V'=a_2+\ldots+a_n-1$, and by induction $R_1$ is at most a simple ramification point of $V'$. Then, the orders of vanishing at $R_1$ of the sections in $V_1$ are 
\begin{center}
$0,\ldots,a_1-1,a_1+1,\ldots,a_1+\ldots+a_n-1$ or

$0,\ldots,a_1-1,a_1+1,\ldots,a_1+\ldots+a_n-2,a_1+\ldots+a_n$,
\end{center}
i.e., $\{0,\ldots,a_1-1,a_1+1,\ldots,a_1+\ldots+a_n\}-\{l\}$, where $l=a_1+\ldots+a_n-1$ or $l=a_1+\ldots+a_n$. Then
\begin{equation*}
wt_{V_1}(R_1)=wt_{V}(R_1)+a_1+\ldots+a_n-1-l=(\sum_{j\neq 1}a_j)-1+\epsilon_1,
\end{equation*}
with $\epsilon_1=a_1+\ldots+a_n-l\in \{0,1\}$. Analogously, we have for each $i$
\begin{center}
$wt_{V_1}(R_i)=(\sum\limits_{j\neq i} a_{j})-1+\epsilon_i$, where $\epsilon_i\in \{0,1\}$.
\end{center}
Finally, using the equalities we have shown, we get
\begin{center}
$\sum wt_{V_1}(R_i)=(a_1+\ldots+a_n-1)(n-1)-1+\sum \epsilon_i$
\end{center}
On the other hand, by Plücker formula, deg$(R_{V_1})=(a_1+\ldots+a_n-1)(n-1)$. Then $0\leq \sum \epsilon_i \leq 1$ and $V_1$ has $1-\sum \epsilon_i$ ramification points on $\mathbb{P}^{1}-\{R_1,\ldots,R_n\}$, counted with their respective weights. This proves the proposition.   \hfill $\Box$

\section{Preliminaries on elliptic curves}\label{curvaeliptica}

\begin{prop}\label{elliptic}
Let $E$ be a smooth elliptic curve, $A$ a point of $E$ and $g\geq 3$ an odd positive integer. Let $\mathcal{L}:=\cal{O}_{E}((2g-2)A)$. Consider the linear system
\begin{center}
$V:=H^{0}(\cal{O}_{E}(gA))\subseteq H^{0}(\mathcal{L})$.
\end{center}
Let $V_{1}\subseteq H^{0}(\mathcal{L})$ be a linear system of dimension $g-1$ such that
\begin{center}
$H^{0}(\cal{O}_{E}((g-2)A))\subseteq V_{1}\subseteq V$.
\end{center}
Then $wt_{V_{1}}(Q)\leq 2$ for every $Q\in E-\{A\}$ and $wt_{V_1}(A)=(g-1)^2+\epsilon$, where $\epsilon\in \{0,1\}$.
\end{prop}
{\em Proof.} Let $Q\in E-\{A\}$. Notice that, since $E$ is elliptic, we have that dim$_{\mathbb{C}}V=g$, $V(-gA)=H^{0}(\cal{O}_{E}((g-2)A))$ has dimension $g-2$, $V(-(g-3)Q)=H^{0}(\cal{O}_{E}(gA-(g-3)Q))$ has dimension $3$ and $V(-gA)\cap V(-(g-3)Q)=H^{0}(\cal{O}_{E}((g-2)A-(g-3)Q))$ has dimension $1$. Then, by dimension considerations
\begin{center}
$V=V(-gA)+V(-(g-3)Q)$. 
\end{center}
Since by hypothesis $V(-gA)$ is contained in $V_{1}$, we get $V(-(g-3)Q)\nsubseteq V_{1}$ and hence dim$_{\mathbb{C}}V_{1}(-(g-3)Q)=2$. On the other hand, the orders of vanishing at $Q$ of the sections in $V$ are of the form
$0,\ldots,g-2,a_{g-1}$, where $g-1\leq a_{g-1}\leq g$, and hence the orders of vanishing at $Q$ of the sections in $V_{1}$ are of the form $\{0,\ldots,g-2,a_{g-1}\}-\{l\}$, for some $l$. As dim$_{\mathbb{C}}V_{1}(-(g-3)Q)=2$, it follows that $l\geq g-3$. Thus $wt_{V_{1}}(Q)=wt_{V}(Q)+g-1-l\leq 2$, if $Q$ is an ordinary point of $V$.

Now, assume $Q$ is an ordinary point of $V(-gA)$. Then, by dimension considerations
\begin{center}
$V=V(-gA)\oplus V(-(g-2)Q)$.
\end{center}
Thus $V(-(g-2)Q)\nsubseteq V_{1}$, and hence dim$_{\mathbb{C}}V_{1}(-(g-2)Q)=1$. It follows that $l\geq g-2$, and since $wt_{V}(Q)\leq 1$, we get $wt_{V_{1}}(Q)=wt_{V}(Q)+g-1-l\leq 2$.

Now, we will prove that $H^{0}(\cal{O}_{E}(gA))$ and $H^{0}(\cal{O}_{E}((g-2)A))$ do not have ramification points in common on $E-\{A\}$. Suppose by contradiction that there exists $Q\in E-\{A\}$ which is a ramification point in common of both $H^{0}(\cal{O}_{E}(gA))$ and $H^{0}(\cal{O}_{E}((g-2)A))$. Then $h^{0}(\cal{O}_{E}(gA-gQ))=1$ and $h^{0}(\cal{O}_{E}((g-2)A-(g-2)Q))=1$. Thus $gA$ and $gQ$ are linearly equivalent divisors and the same property is true for $(g-2)A$ and $(g-2)Q$. Therefore, $2A$ and $2Q$ are linearly equivalent divisors. Now let $g=2n+1$; since $2A$ and $2Q$ are linearly equivalent divisors, we have that $2nA$ and $2nQ$ are linearly equivalent divisors. As $(2n+1)A$ and $(2n+1)Q$ are linearly equivalent divisors, it follows that $A$ and $Q$ are linearly equivalent divisors and hence $Q=A$, a contradiction. 

Finally, we will compute $wt_{V_1}(A)$. Since the orders of vanishing at $A$ of the sections in $V$ are $g-2,\ldots,2g-4,2g-2$,
we have that the orders of vanishing at $A$ of the sections in $V_1$ are of the form
\begin{center}
$\{g-2,\ldots,2g-4,2g-2\}-\{l\}$
\end{center}
for some $l$. Since $V(-gA)\subseteq V_1\subseteq V$, we have that $V_1(-gA)=V(-gA)$. Then dim$_{\mathbb{C}}V_1(-gA)=g-2$, and hence $g-2\leq l\leq g-1$. So, since $wt_{V}(A)=(g-1)^2$,
\begin{center}
$wt_{V_1}(A)=wt_{V}(A)+g-1-l=(g-1)^2+\epsilon$, with $\epsilon=g-1-l\in \{0,1\}$.
\end{center}
  \hfill $\Box$

\section{The divisor}\label{principal}

%\section{Introduction}
\begin{subs}{\em(}The divisor{\em)}\label{coeffa} {\em
Our aim is to study the class of the divisor $\overline{S^{2}W}$ in Pic$_{\text{fun}}(\overline{M_{g}})$, defined as
the closure of the locus of smooth curves $C$ with a pair of points $(P,Q)$ satisfying that $Q$ is a ramification point of the linear system $H^{0}(\omega_{C}(-P))$ with ramification weight at least $3$. 

Write the class of the divisor as 
\begin{center}
$\overline{S^{2}W}=a \lambda-a_{0}\delta_{0}-a_{1}\delta_{1}-\ldots-a_{[g/2]}\delta_{[g/2]}$.
\end{center}
First, we will compute the coefficient $a$.
Let $\pi:\mathcal{X}\rightarrow T$ be a family of smooth curves over a smooth curve $T$. Consider the double product $\mathcal{Y}:=\mathcal{X}\times_{T}\mathcal{X}$
as a family of curves via the first projection $p_{1}:\mathcal{Y}\rightarrow \mathcal{X}$.
Let $W$ be the ramification divisor of the invertible sheaf $\mathcal{L}:=\omega_{p_{1}}(-\Delta)$ with respect to $p_{1}$.
Notice that $h^{0}(\mathcal{L}\big|_{\mathcal{Y}_{P}})=h^{0}(\omega_{\mathcal{X}_{\pi(P)}}(-P))=g-1$ for every $P\in \mathcal{X}$. Then $p_{1*}(\mathcal{L})$ is locally free of rank $g-1$.

Now, we will compute $\pi_{*}p_{1*}([S^{2}W])$. By the Thom-Porteous formula:
\begin{equation}\label{divisor1}
[W]=c_{1}(J_{p_{1}}^{g-2}(\mathcal{L}))-c_{1}(p_{1}^{*}p_{1*}(\mathcal{L})).
\end{equation}
By using the truncation exact sequences and the Whitney formula, we obtain
\begin{equation}\label{divisor2}
c_{1}(J_{p_{1}}^{g-2}(\mathcal{L}))={g-1 \choose 2}c_{1}(\omega_{p_{1}})+(g-1)c_{1}(\mathcal{L}).
\end{equation}
We have to compute $c_{1}(p_{1*}(\mathcal{L}))$.
Notice that by Riemann-Roch we have $h^{1}(\mathcal{L}\big|_{\mathcal{Y}_{P}})=1$ for every $P\in \mathcal{X}$, as $h^{0}(\mathcal{L}\big|_{\mathcal{Y}_{P}})=g-1$.
It follows that $R^{1}p_{1*}(\mathcal{L})$ is invertible.

Consider the long exact sequence
\begin{center}
$0\rightarrow p_{1*}(\mathcal{L})\rightarrow p_{1*}(\omega_{p_{1}})\rightarrow p_{1*}(\omega_{p_{1}}\big|_{\Delta})\rightarrow $
\end{center}

\begin{center}
$R^{1}p_{1*}(\mathcal{L})\rightarrow R^{1}p_{1*}(\omega_{p_{1}})\rightarrow R^{1}p_{1*}(\omega_{p_{1}}\big|_{\Delta})\rightarrow 0$.
\end{center}
Since $R^{1}p_{1*}(\omega_{p_{1}}\big|_{\Delta})=0$, as the restriction of $\omega_{p_{1}}\big|_{\Delta}$ to each fiber is supported at a point, we have a surjection
$R^{1}p_{1*}(\mathcal{L})\rightarrow R^{1}p_{1*}(\omega_{p_{1}})$. As $R^{1}p_{1*}(\mathcal{L})$ is an invertible sheaf and $R^{1}p_{1*}(\omega_{p_{1}})\cong \cal{O}_{\mathcal{X}}$, it follows that $R^{1}p_{1*}(\mathcal{L})\cong R^{1}p_{1*}(\omega_{p_{1}})$. Then we have an exact sequence
\begin{center}
$0\rightarrow p_{1*}(\mathcal{L})\rightarrow p_{1*}(\omega_{p_{1}})\rightarrow p_{1*}(\omega_{p_{1}}\big|_{\Delta})\rightarrow 0$.
\end{center}
%Por Grothendieck-Riemann-Roch temos:

%$ch(p_{1!}(\mathcal{L}))=p_{1*}(ch(\mathcal{L})\cdot td(\mathcal{T}_{\mathcal{Y}/\mathcal{X}}))$

%\begin{center}
%$=p_{1*}((1+c_{1}(\mathcal{L})+\frac{c_{1}(\mathcal{L})^{2}}{2}+...)\cdot (1-\frac{K_{p_{1}}}{2}+td_{2}(\mathcal{T}_{\mathcal{Y}/\mathcal{X}})+...))$
%\end{center}

%\begin{center}
%$=p_{1*}(1+(c_{1}(\mathcal{L})-\frac{K_{p_{1}}}{2})+(\frac{c_{1}(\mathcal{L})^{2}}{2}-\frac{K_{p_{1}}c_{1}(\mathcal{L})}{2}+td_{2}(\mathcal{T}_{\mathcal{Y}/\mathcal{X}}))+...)$
%\end{center}

%\begin{center}
%$=p_{1*}(1+(K_{p_{1}}-\Delta-\frac{K_{p_{1}}}{2})+(\frac{K_{p_{1}}^{2}-2K_{p_{1}}\cdot \Delta+\Delta^{2}}{2}-\frac{K_{p_{1}}\cdot(K_{p_{1}}-\Delta)}{2}+td_{2}(\mathcal{T}_{\mathcal{Y}/\mathcal{X}}))+...)$
%\end{center}

%\begin{center}
%$=p_{1*}(1+(\frac{K_{p_{1}}}{2}-\Delta)+(\frac{-K_{p_{1}}\cdot \Delta+\Delta^{2}}{2}+td_{2}(\mathcal{T}_{\mathcal{Y}/\mathcal{X}}))+...)$
%\end{center}
%Lembrando que $p_{1*}(td_{2}(\mathcal{T}_{\mathcal{Y}/\mathcal{X}}))=\lambda_{p_{1}}=\pi^{*}\lambda$, temos
%\begin{center}
%$ch(p_{1!}(\mathcal{L}))=(g-2)+(-K_{\pi}+\pi^{*}\lambda)+...$
%\end{center}
Via the Whitney formula, we have
\begin{center}
$c_{1}(p_{1*}(\mathcal{L}))=c_{1}(p_{1*}(\omega_{p_{1}}))-c_{1}(p_{1*}(\omega_{p_{1}}\big|_{\Delta}))$.
\end{center}
From $\omega_{p_{1}}=p_{2}^{*}\omega_{\pi}$, we get $p_{1*}(\omega_{p_{1}}\big|_{\Delta})=\omega_{\pi}$, and since $p_{1*}(\omega_{p_{1}})=p_{1*}(p_{2}^{*}(\omega_{\pi}))\cong \pi^{*}\pi_{*}\omega_{\pi}$,
\begin{center}
$c_{1}(p_{1*}(\omega_{p_{1}}))=\pi^{*}c_{1}(\pi_{*}\omega_{\pi})=\pi^{*}c_{1}($det $\pi_{*}\omega_{\pi})=\pi^{*}\lambda_{\pi}$.
\end{center}
Therefore
\begin{equation}\label{divisor3}
c_{1}(p_{1*}(\mathcal{L}))=\pi^{*}\lambda-K_{\pi}, 
\end{equation}
where $\lambda:=\lambda_{\pi}$ and $K_{\pi}:=c_{1}(\omega_{\pi})$. Now, let $K_{p_{1}}:=p_{2}^{*}K_{\pi}$ and $K_{p_{2}}:=p_{1}^{*}K_{\pi}$. Then, by (\ref{divisor1}), (\ref{divisor2}) and (\ref{divisor3}),
\begin{equation*}
[W]={g-1 \choose 2}K_{p_{1}}+(g-1)c_{1}(\mathcal{L})-p_{1}^{*}(\pi^{*}\lambda-K_{\pi})={g \choose 2}K_{p_{1}}+K_{p_{2}}-(g-1)\Delta-p_{1}^{*}\pi^{*}\lambda.
\end{equation*}
By the Thom-Porteous formula:
\begin{center}
$[S^{2}W]=c_{3}(J_{p_{1}}^{2}(\cal{O}_{\mathcal{Y}}(W)))$.
\end{center}
Using the truncation exact sequence
\[0\rightarrow \omega_{p_{1}}^{\otimes 2}\otimes \cal{O}_{\mathcal{Y}}(W)\rightarrow J_{p_{1}}^{2}(\cal{O}_{\mathcal{Y}}(W))\rightarrow J_{p_{1}}^{1}(\cal{O}_{\mathcal{Y}}(W))\rightarrow 0\]
and recalling that $J_{p_{1}}^{1}(\cal{O}_{\mathcal{Y}}(W))$ is locally free of rank $2$, we get
\begin{center}
$c_{3}(J_{p_{1}}^{2}(\cal{O}_{\mathcal{Y}}(W)))=c_{2}(J_{p_{1}}^{1}(\cal{O}_{\mathcal{Y}}(W)))c_{1}(\omega_{p_{1}}^{\otimes 2}\otimes \cal{O}_{\mathcal{Y}}(W))$,
\end{center}
and using the truncation exact sequence
\[0\rightarrow \omega_{p_{1}}\otimes \cal{O}_{\mathcal{Y}}(W)\rightarrow J_{p_{1}}^{1}(\cal{O}_{\mathcal{Y}}(W))\rightarrow \cal{O}_{\mathcal{Y}}(W)\rightarrow 0\]
we get $c_{2}(J_{p_{1}}^{1}(\cal{O}_{\mathcal{Y}}(W)))=c_{1}(\cal{O}_{\mathcal{Y}}(W))c_{1}(\omega_{p_{1}}\otimes \cal{O}_{\mathcal{Y}}(W))$.

Therefore
\begin{center}
$[S^{2}W]=c_{3}(J_{p_{1}}^{2}(\cal{O}_{\mathcal{Y}}(W)))=[W](K_{p_{1}}+[W])(2K_{p_{1}}+[W])$.
\end{center}

On the other hand, since $\cal{O}(-\Delta)\big|_{\Delta}\cong \omega_{\pi}$ (identifying $\Delta$ with $\mathcal{X}$), we have that $\cal{O}(-\Delta)\big|_{\Delta}=(p_{2}^{*}\omega_{\pi})\big|_{\Delta}=\omega_{p_{1}}\big|_{\Delta}$ and $\cal{O}(-\Delta)\big|_{\Delta}=(p_{1}^{*}\omega_{\pi})\big|_{\Delta}=\omega_{p_{2}}\big|_{\Delta}$. It follows that $\Delta^{2}=-K_{p_{1}}\cdot \Delta=-K_{p_{2}}\cdot \Delta$.
Using the projection formula and the following formulas

(1)$K_{p_{1}}^{3}=0$, $K_{p_{2}}^{3}=0$ and $(p_{1}^{*}\pi^{*}\lambda)^{2}=0$.

(2)$p_{1*}(K_{p_{1}}\cdot \Delta)=K_{\pi}$, $\Delta^{2}=-K_{p_{1}}\cdot \Delta=-K_{p_{2}}\cdot \Delta$.

(3)$\pi_{*}(K_{\pi}^{2})=12\lambda$, $\pi_{*}(K_{\pi})=2g-2$.

(4)$\pi^{*}\pi_{*}(\alpha)=p_{1*}p_{2}^{*}(\alpha)$ for every cycle $\alpha$ on $\mathcal{X}$, we get:
\begin{center}
$\pi_{*}p_{1*}([S^{2}W])=(9g^5-51g^4+129g^3-207g^2+174g-54)\lambda$.
\end{center}

Therefore, 
\begin{center}
$a=9g^5-51g^4+129g^3-207g^2+174g-54$. 
\end{center}

Now, as it will be important for us, we will recall the definition of a flag curve. A \textit{flag curve} is a nodal curve $X$ satisfying the following properties:

$(1)$ It is of \textit{compact type}, i.e., the number of nodes of $X$ is smaller (by one) than the number of components.

$(2)$ Each component of $X$ is either $\mathbb{P}^{1}$ or an elliptic curve.

$(3)$ Each elliptic component of $X$ contains exactly one node of $X$.

$(4)$ Each $\mathbb{P}^{1}$ contains at least $2$ nodes of $X$.
}
\end{subs}

%For $g$ odd and $g\geq 5$, we will obtain the coefficient $a_{i}$ in terms of the coefficient $a_{1}$ for every $i>1$, in Chapter $5$, by using the method of test curves. We will use $[g/2]-1$ test curves, which are induced by families of flag stable curves over $\mathbb{P}^{1}$. Of crucial importance in the use of the test curves is Proposition $5.2.1$, which is similar to $[HMo]$, Thm $6.65$, $(2)$; in fact, Proposition $5.2.1$ implies that result. To be able to apply our Proposition $5.2.1$, we will use Proposition $5.1.1$, which is a general result about flag curves. We end up with (see Chapter $5$ for more details)
%\begin{center}
%$a_i=(i(g-i)/(g-1))a_1$, for each $2\leq i\leq [g/2]$.
%\end{center}

%Our main theorem is:

\begin{prop}\label{propflag}
Let $X$ be a flag curve of genus $g$. Assume $g\geq 3$ is an odd integer and let $p:\cal{C}\rightarrow \Sigma:=Spec\,\mathbb{C}[[t]]$ be a smoothing of $X$. Let $\cal{C}_{*}$ be the generic fiber of $p$ and $\overline{\cal{C}}_{*}$ the geometric generic fiber. Then $\overline{\cal{C}}_{*}$ satisfies the following condition:

for each $P_{*}\in \overline{\cal{C}}_{*}$, the ramification points of the complete linear system $H^{0}(\omega_{\overline{\cal{C}}_{*}}(-P_{*}))$ have ramification weight at most $2$.
\end{prop}
{\em Proof.} Let $P_{*}\in \overline{\cal{C}}_{*}$. After base change, we may assume that $P_{*}$ is a rational point of $\cal{C}_{*}$, and thus there is a section $\Gamma$ of $p$ intersecting $\cal{C}_{*}$ at $P_{*}$. After base changes and a sequence of blowups at the singular points of the special fiber $\cal{C}_{0}$, we may assume that $\cal{C}$ is regular and that $\Gamma$ intersects the special fiber at a point $P$ which is not a node of $\cal{C}_{0}$. After all the base changes and the sequence of blowups, each node is replaced by a chain of rational smooth curves and $\cal{C}_{0}$ is still a flag curve.

Let $\mathcal{L}:=\omega_{p}(-\Gamma)$ and let $W_0$ be the limit ramification divisor of $\mathcal{L}$. To prove the statement of the proposition, it is enough to show that mult$_{Q}(W_0)\leq 2$ for every $Q\in \cal{C}_{0}$. There are two cases to consider.
%the limit ramification divisor $W_0$ of $\mathcal{L}$ has multiplicity at most $2$ at each point $Q\in \cal{C}_{0}$. 

Case $(1)$: $P$ lies on a rational component $Y$ of $\cal{C}_{0}$. 

We will show that mult$_{Q}(W_0)\leq 2$ for every $Q\in \cal{C}_{0}$. To prove this, we will show that the multiplicity of $W_0$ at each node of $\cal{C}_{0}$ is $0$, and mult$_{Q}(W_0)\leq 2$ if $Q\in \cal{C}_{0}$ is not a node. 

The limit linear system of $\omega_{p}$ on $Y$ is of the form (see \cite{EH2})
\begin{center}
$V:=H^{0}(\omega_{Y}((a_{1}+1)R_1))\oplus \ldots \oplus H^{0}(\omega_{Y}((a_{n}+1)R_n))$

$\subseteq H^{0}(\omega_{Y}(2a_1R_{1}+\ldots +2a_nR_n))$,
\end{center}
where $n$ is the number of connected components of $\cal{C}_{0}-Y$, the integers $a_j$ are the genera of the closures of the connected components of $\cal{C}_{0}-Y$, and each $R_j$ is the point of intersection of $Y$ and the connected component of the closure of $\cal{C}_{0}-Y$ of genus $a_j$.
Notice that if $\omega_{p}(D_Z)$ has degree $2g-2$ on a component $Z$ of $\cal{C}_{0}$ and degree $0$ on the other components of $\cal{C}_{0}$, where $D_Z\subseteq \cal{C}_{0}$ is a divisor, then $\mathcal{L}(D_Z)$ has focus on $Z$. In this way, we can get a limit linear system $V_Z$ of $\mathcal{L}$ on each component $Z$ of $\cal{C}_{0}$, and the connecting number between $\mathcal{L}(D_{Z_1})$ and $\mathcal{L}(D_{Z_2})$ corresponding to components $Z_1\neq Z_2$ of $\cal{C}_{0}$ is equal to the connecting number between $\omega_{p}(D_{Z_1})$ and $\omega_{p}(D_{Z_2})$ corresponding to $Z_1$ and $Z_2$. 
The limit linear system of $\mathcal{L}$ on $Y$ is
\begin{center}
$V_Y=V(-P)\subseteq H^{0}(\omega_{Y}(-P+2a_1R_{1}+\ldots +2a_nR_n))$.
\end{center}
It follows from Proposition \ref{rationalp} that $wt_{V_Y}(Q)\leq 1$ if $Q\in Y-\{R_1,\ldots, R_n\}$. In particular, mult$_{Q}(W_0)\leq 2$ if $Q\in Y$ is not a node of $\cal{C}_{0}$. Now, we will prove that the multiplicity of $W_0$ at each point $R_j$ is $0$.  We may assume $j=1$. Now, assume that $R_1$ is the point of intersection of $Y$ and a rational component $Y_1$ of $\cal{C}_{0}$. Since $a_1+\ldots+a_n=g$ and $V_Y$ has dimension $g-1$, it follows from Proposition \ref{rationalp} that
\begin{equation*}
wt_{V_Y}(R_1)=g-a_1-1+\epsilon_1+(a_1-1)(g-1)
\end{equation*}
where $\epsilon_1\in \{0,1\}$. Let $R'_1,\ldots,R'_m$ be the nodes of $\cal{C}_{0}$ lying on $Y_1$.  We may assume $R'_m=R_1$. Let $V_{Y_1}$ be the limit linear system of $\mathcal{L}$ on $Y_1$. By Formula (\ref{limitramif}), we have
\begin{align}
\text{mult}_{R_1}(W_0)&=wt_{V_Y}(R_1)+wt_{V_{Y_1}}(R'_m)+(g-1)(g-2-(2g-2))\nonumber \\
&=wt_{V_{Y_1}}(R'_m)+(g-1)(a_1-g)-a_1+\epsilon_1. \label{flag1}
\end{align}
On the other hand, the limit linear system of $\omega_{p}$ on $Y_1$ is of the form
\begin{center}
$V':=H^{0}(\omega_{Y_1}((a'_{1}+1)R'_1))\oplus \ldots \oplus H^{0}(\omega_{Y_1}((a'_{m}+1)R'_m))$

$\subseteq H^{0}(\omega_{Y_1}(2a'_1R'_{1}+\ldots +2a'_mR'_m))$.
\end{center}
We have that $V_{Y_1}\subseteq V'$. It follows from Proposition \ref{rational} that the orders of vanishing at $R'_m$ of the sections in $V'$ are 
\begin{center}
$0+(a'_m-1),\ldots,a'_m-1+(a'_m-1),a'_m+1+(a'_m-1),\ldots,a'_1+\ldots+a'_m+(a'_m-1)$
\end{center}
and $wt_{V'}(R'_m)=(\sum_{j\neq m}a'_j)+(a'_m-1)g$. Thus, the orders of vanishing at $R'_m$ of the sections in $V_{Y_1}$ are of the form
\begin{center}
$\{a'_m-1,\ldots,2(a'_m-1),a'_m+1+(a'_m-1),\ldots,a'_1+\ldots+a'_m+(a'_m-1)\}-\{l\}$,
\end{center}
for some $l$. Thus 
\begin{equation}\label{flag2}
wt_{V_{Y_1}}(R'_m)=wt_{V'}(R'_m)+g-1-l=g-a'_m+(a'_m-1)g+g-1-l.
\end{equation}
Then, since $a_1+a'_m=g$, by (\ref{flag1}) and (\ref{flag2}) we have
\begin{equation*}
\begin{split}
\text{mult}_{R_1}(W_0)&=g-a'_m+(a'_m-1)g+g-1-l+(g-1)(a_1-g)-a_1+\epsilon_1\\
&=(a'_m-1-l)+\epsilon_1\leq \epsilon_1\leq 1.
\end{split}
\end{equation*}
Since the intersection multiplicity of the ramification divisor of $\mathcal{L}$ and the special fiber at the node $R_1$ cannot be $1$, we have mult$_{R_1}(W_0)=0$.
(Notice that, the only important information about $V_Y$ we have used in the reasoning above is the ramification weight of $V_Y$ at the point $R_1$.)

Now, we are going to see what happens on $Y_1$. We have to prove that the multiplicity of $W_0$ at each point $R'_j$ is $0$, and mult$_{Q}(W_0)\leq 2$ if $Q\in Y_1$ is not a node of $\cal{C}_{0}$. Since mult$_{R_1}(W_0)=0$, $a'_m-1-l=-\epsilon_1$. This implies that $l=a'_m-1$ or $l=a'_m$. Then dim$_{\mathbb{C}}V_{Y_1}(-(a'_m+1)R'_m)=g-2$ and hence $V_{Y_1}\supseteq V'(-(a'_m+1)R'_m)$, i.e.,
\begin{center}
$V_{Y_1}\supseteq H^{0}(\omega_{Y_1}((a'_{1}+1)R'_1))\oplus \ldots \oplus H^{0}(\omega_{Y_1}((a'_{m}-1)R'_m))$.
\end{center}
Now, using Proposition \ref{rationalr} for $V_1$ equal to $V_{Y_1}$, we are able to use the same reasoning above to conclude that mult$_{R'_j}(W_0)=0$ if $R'_j$ is the point of intersection of $Y_1$ and a rational component of $\cal{C}_{0}$. Also, using Proposition \ref{rationalr}, we get that mult$_{Q}(W_0)\leq 2$ if $Q\in Y_1$ is not a node of $\cal{C}_{0}$. Notice that, we can use the same reasoning above, repeatedly, for each rational component in $\cal{C}_{0}$.

It remains to prove that, if $\bar{Y}$ is a rational component of $\cal{C}_{0}$ intersecting an elliptic component $E$ of $\cal{C}_{0}$, then the point of intersection of $\bar{Y}$ and $E$ does not appear in $W_0$ and mult$_{Q}(W_0)\leq 2$ if $Q\in E$ is not a node of $\cal{C}_{0}$. Let $\bar{R}_1,\ldots,\bar{R}_k$ be the nodes of $\cal{C}_0$ lying on $\bar{Y}$ and $A$ the node of $\cal{C}_0$ lying on $E$. We may assume $\bar{R}_k=A$. Let $V_{\bar{Y}},V_E$ be the limit linear systems of $\mathcal{L}$ on $\bar{Y}$ and $E$ respectively. We have an equality of the form
\begin{center}
$wt_{V_{\bar{Y}}}(\bar{R}_k)=(\sum\limits_{j\neq k}\bar{a}_j)-1+\epsilon_k+(\bar{a}_k-1)(g-1)$,
\end{center}
where $\epsilon_k\in \{0,1\}$, the integers $\bar{a}_j$ are the genera of the closures of the connected components of $\cal{C}_{0}-\bar{Y}$, and each $\bar{R}_j$ is the point of intersection of $\bar{Y}$ and the connected component of the closure of $\cal{C}_{0}-\bar{Y}$ of genus $\bar{a}_j$. Since $\bar{a}_k=1$, $wt_{V_{\bar{Y}}}(\bar{R}_k)=g-2+\epsilon_k$. Then, by Formula (\ref{limitramif}), we have
\begin{align}
\text{mult}_{\bar{R}_k}(W_0)&=wt_{V_{\bar{Y}}}(\bar{R}_k)+wt_{V_{E}}(A)+(g-1)(g-2-(2g-2))\nonumber \\
&=wt_{V_{E}}(A)+g-2-g(g-1)+\epsilon_k. \label{flag3}
\end{align}
On the other hand, the limit linear system of $\omega_{p}$ on $E$ is
\begin{center}
$V':=H^{0}(\omega_{E}(gA))\subseteq H^{0}(\omega_{E}((2g-2)A))$.
\end{center}
It follows that $V_{E}\subseteq H^{0}(\omega_{E}(gA))$. The orders of vanishing at $A$ of the sections in $V'$ are 
\begin{center}
$g-2,\ldots,2g-4,2g-2$
\end{center}
and $wt_{V'}(A)=(g-1)^2$. Thus, the orders of vanishing at $A$ of the sections in $V_{E}$ are of the form 
\begin{center}
$\{g-2,\ldots,2g-4,2g-2\}-\{l\}$, for some $l$.
\end{center}
Then
\begin{equation}\label{flag4}
wt_{V_{E}}(A)=wt_{V'}(A)+g-1-l=g^2-g-l.
\end{equation}
By (\ref{flag3}) and (\ref{flag4}), we get
\begin{equation*}
\text{mult}_{\bar{R}_k}(W_0)=g^2-g-l+g-2-g(g-1)+\epsilon_k=(g-2-l)+\epsilon_k\leq \epsilon_k\leq 1.
\end{equation*}
It follows that mult$_{\bar{R}_k}(W_0)=0$ and hence $\bar{R}_k$ is not a limit ramification point, and $l=g-2$ or $l=g-1$, which implies that dim$_{\mathbb{C}}V_E(-gA)=g-2$ and hence $V_E$ contains $V'(-gA)=H^{0}(\omega_{E}((g-2)A))$. By using Proposition \ref{elliptic}, we get mult$_{Q}(W_0)\leq 2$ if $Q\in E$ is not a node of $\cal{C}_{0}$. This proves the case $(1)$.

Case $(2)$: $P$ lies on an elliptic component $E$ of $\cal{C}_{0}$.

We will show that mult$_{Q}(W_0)\leq 2$ for every $Q\in \cal{C}_{0}$. To prove this, we will show that the multiplicity of $W_0$ at each node of $\cal{C}_{0}$ is $0$, and mult$_{Q}(W_0)\leq 2$ if $Q\in \cal{C}_{0}$ is not a node. 

Let $A$ be the node of $\cal{C}_0$ lying on $E$. Since the limit linear system of $\omega_p$ on $E$ is
\begin{center}
$V:=H^{0}(\omega_{E}(gA))\subseteq H^{0}(\omega_{E}((2g-2)A))$,
\end{center}
the limit linear system of $\mathcal{L}$ on $E$ is
\begin{center}
$V_E=V(-P)=H^{0}(\omega_{E}(gA-P))\subseteq H^{0}(\omega_{E}((2g-2)A-P))$.
\end{center}
Notice that $V_E$ has at most simple ramification points on $E-\{A\}$, whence mult$_{Q}(W_0)\leq 2$ if $Q\in E$ is not a node of $\cal{C}_{0}$. Now, we have to prove that the multiplicity of $W_0$ at the point $A$ is $0$. We have $V_E(-nA)=V_E$ for $0\leq n\leq g-2$. Also, for $n\geq g-2$, $V_E(-nA)=H^{0}(\omega_{E}((2g-2)A-P-nA))$. Then 
\begin{center}
dim$_{\mathbb{C}}V_E(-nA)=2g-3-n$ for every $g-2\leq n\leq 2g-4$.
\end{center} 
Thus, since $V_E(-(2g-3)A)=H^{0}(\omega_E(A-P))=0$, the orders of vanishing at $A$ of the sections in $V_E$ are $g-2,\ldots,2g-4$, and hence 
\begin{equation}\label{flag5}
wt_{V_E}(A)=(g-1)(g-2).
\end{equation}
Let $Y$ be the rational component of $\cal{C}_{0}$ intersecting $E$ at the point $A$. Let $R_1,\ldots,R_n$ be the nodes of $\cal{C}_0$ lying on $Y$. We may assume $R_n=A$. Let $V_Y$ be the limit linear system of $\mathcal{L}$ on $Y$. The limit linear system of $\omega_p$ on $Y$ is of the form
\begin{center}
$V':=H^{0}(\omega_{Y}((a_{1}+1)R_1))\oplus \ldots \oplus H^{0}(\omega_{Y}((a_{n}+1)R_n))$

$\subseteq H^{0}(\omega_{Y}(2a_1R_{1}+\ldots +2a_nR_n))$,
\end{center}
and we have $V_Y\subseteq V'$. It follows from Proposition \ref{rational} that the orders of vanishing at $R_n$ of the sections in $V'$ are 
\begin{center}
$0+(a_n-1),\ldots,a_n-1+(a_n-1),a_n+1+(a_n-1),\ldots,a_1+\ldots+a_n+(a_n-1)$
\end{center}
and $wt_{V'}(R_n)=(\sum\limits_{j\neq n}a_j)+(a_n-1)g$. Thus, the orders of vanishing at $R_n$ of the sections in $V_{Y}$ are of the form
\begin{center}
$\{a_n-1,\ldots,2(a_n-1),a_n+1+(a_n-1),\ldots,a_1+\ldots+a_n+(a_n-1)\}-\{l\}$,
\end{center}
for some $l$. Then, since $a_n=1$,
\begin{equation}\label{flag6}
wt_{V_{Y}}(R_n)=wt_{V'}(R_n)+g-1-l=g-a_n+(a_n-1)g+g-1-l=2g-2-l,
\end{equation}
By Formula (\ref{limitramif}), (\ref{flag5}) and (\ref{flag6}), we get
\begin{equation*}
\text{mult}_{C}(W_0)=wt_{V_{Y}}(R_n)+wt_{V_E}(A)+g(g-1-(2g-2))=-l.
\end{equation*}
It follows that mult$_{C}(W_0)=0$ and hence $A$ is not a limit ramification point, and $l=0$, which implies that 
\begin{center}
$V_Y=H^{0}(\omega_{Y}((a_{1}+1)R_1))\oplus \ldots \oplus H^{0}(\omega_{Y}((a_{n-1}+1)R_{n-1}))$.
\end{center}
It follows from Proposition \ref{rational} that, for every $k\neq n$
\begin{equation*}
wt_{V_Y}(R_k)=(\sum_{j\neq k,n}a_j)+(a_k-1)(g-1)=(\sum_{j\neq k}a_j)-1+(a_k-1)(g-1).
\end{equation*}
Thus, the proof of this case follows as in the case $(1)$.   \hfill $\Box$

\begin{subs}{\em(}Families of flag curves over $\mathbb{P}^{1}${\em)} {\em Let $g$ and $i$ be positive integers such that $g\geq 5$ and $2\leq i\leq [g/2]$. We will define a family of curves over $\mathbb{P}^{1}$ in the following steps:

{\em Step }$1$: Fix $3$ distinct points $R,S$ and $T$ on $\mathbb{P}^{1}$, and let $\mathbb{P}^{1}_{i-1}:=\mathbb{P}^{1}$, $R_{i-1}:=R$, $S_{i-1}:=S$ and $T_{i-1}:=T$ (we use this notation to extend a notation we will see later). Begin with the fibered product $\mathbb{P}^{1}_{i-1}\times \mathbb{P}^{1}_{i-1}$, and then blow up the points $(R_{i-1},R_{i-1}),(S_{i-1},S_{i-1})$ and $(T_{i-1},T_{i-1})$. Let $\mathbb{P}^{1}_{(R_{i-1},R_{i-1})}, \mathbb{P}^{1}_{(S_{i-1},S_{i-1})}$ and $\mathbb{P}^{1}_{(T_{i-1},T_{i-1})}$ be the rational curves on the blowup $(\mathbb{P}^{1}_{i-1}\times \mathbb{P}^{1}_{i-1})^{\widetilde{}}$ over the points $(R_{i-1},R_{i-1}),(S_{i-1},S_{i-1}), (T_{i-1},T_{i-1})$ of $\mathbb{P}^{1}_{i-1}\times \mathbb{P}^{1}_{i-1}$. The points in the intersections $(\{R_{i-1}\}\times \mathbb{P}^{1}_{i-1})^{\widetilde{}}\cap \mathbb{P}^{1}_{(R_{i-1},R_{i-1})}$, $(\{S_{i-1}\}\times \mathbb{P}^{1}_{i-1})^{\widetilde{}}\cap \mathbb{P}^{1}_{(S_{i-1},S_{i-1})}$ and $(\{T_{i-1}\}\times \mathbb{P}^{1}_{i-1})^{\widetilde{}}\cap \mathbb{P}^{1}_{(T_{i-1},T_{i-1})}$ will be denoted $R',S'$ and $T'$ respectively. Also, abusing notation, the strict transform of each fiber $\{Q\}\times \mathbb{P}^{1}_{i-1}\subseteq \mathbb{P}^{1}_{i-1}\times \mathbb{P}^{1}_{i-1}$ will be denoted $\mathbb{P}^{1}_{i-1}$.

{\em Step }$2$: Now, fix $g$ smooth pointed elliptic curves $(E_1,A_1),\ldots,(E_g,A_g)$. Let $\cal{Y}$ be the disjoint union of $(\mathbb{P}^{1}_{i-1}\times \mathbb{P}^{1}_{i-1})^{\widetilde{}}$, $\mathbb{P}^{1}\times E_i$  and $\mathbb{P}^{1}\times E_{i+1}$ modulo the identification of the strict transform of the diagonal $\widetilde{\Delta}\subseteq (\mathbb{P}^{1}_{i-1}\times \mathbb{P}^{1}_{i-1})^{\widetilde{}}$ with $\mathbb{P}^{1}\times \{A_i\}\subseteq \mathbb{P}^{1}\times E_i$, and the identification of the strict transform $(\mathbb{P}^{1}_{i-1}\times \{S_{i-1}\})^{\widetilde{}}\subseteq (\mathbb{P}^{1}_{i-1}\times \mathbb{P}^{1}_{i-1})^{\widetilde{}}$ with $\mathbb{P}^{1}\times \{A_{i+1}\}\subseteq \mathbb{P}^{1}\times E_{i+1}$.

{\em Step }$3$: Assume $i\geq 3$ and consider a chain of $i-2$ three pointed rational curves $(\mathbb{P}^{1}_{1},R_{1},S_{1},T_{1}),\ldots,(\mathbb{P}^{1}_{i-2},R_{i-2},S_{i-2},T_{i-2})$ with $T_j=R_{j+1}$ for every $1\leq j\leq i-3$. Now, attach the elliptic curves $E_1,\ldots,E_{i-1}$ at the points $R_1,S_1,S_2,\ldots,S_{i-2}$ respectively, identifying the points $A_1,\ldots,A_{i-1}$ with the points $R_1,S_1,S_2,\ldots,S_{i-2}$ respectively, obtaining a nodal curve which we will call $X_i$. If $i=2$, we set $X_i:=E_1$ and $T_{i-2}:=A_1$.
Analogously, consider a chain of $g-i-2$ three pointed rational curves $(\mathbb{P}^{1}_{i},R_{i},S_{i},T_{i}),\ldots,(\mathbb{P}^{1}_{g-3},R_{g-3},S_{g-3},T_{g-3})$ such that $T_j=R_{j+1}$ for every $i\leq j\leq g-4$. Now, attach the elliptic curves $E_{i+2},\ldots,E_g$ at the points $S_i,\ldots,S_{g-3},T_{g-3}$ respectively, identifying the points $A_{i+2},\ldots,A_g$ with the points $S_i,\ldots,S_{g-3},T_{g-3}$ respectively, obtaining a nodal curve which we will call $Y_i$. 

{\em Step }$4$: Finally, let $\cal{X}$ be the disjoint union of $\cal{Y}$, $\mathbb{P}^{1}\times X_i$ and $\mathbb{P}^{1}\times Y_i$ modulo the identification of $(\mathbb{P}^{1}_{i-1}\times \{R_{i-1}\})^{\widetilde{}}\subseteq \cal{Y}$ with $\mathbb{P}^{1}\times \{T_{i-2}\}\subseteq \mathbb{P}^{1}\times X_i$, and the identification of $(\mathbb{P}^{1}_{i-1}\times \{T_{i-1}\})^{\widetilde{}}\subseteq \cal{Y}$ with $\mathbb{P}^{1}\times \{R_i\}\subseteq \mathbb{P}^{1}\times Y_i$. This gives a family $\pi_i:\cal{X}\rightarrow \mathbb{P}^{1}$ of stable curves of genus $g$.

Abusing notation, for each fiber $F$ of $\pi_i$, we denote by $R_{i-1},S_{i-1}$ and $T_{i-1}$ the points in the intersections $F\cap (\mathbb{P}^{1}_{i-1}\times \{R_{i-1}\})^{\widetilde{}}$, $F\cap (\mathbb{P}^{1}_{i-1}\times \{S_{i-1}\})^{\widetilde{}}$ and $F\cap (\mathbb{P}^{1}_{i-1}\times \{T_{i-1}\})^{\widetilde{}}$, respectively. 
Figure $1$ describes the family given by $\pi_i$. We denote by $[\pi_i]:\mathbb{P}^{1}\rightarrow \overline{M_{g}}$ the map induced by the family $\pi_i:\cal{X}\rightarrow \mathbb{P}^{1}$.

%\vspace{2cm}

\begin{figure}[h]
	\centering
		\includegraphics[height=13cm,width=16cm]{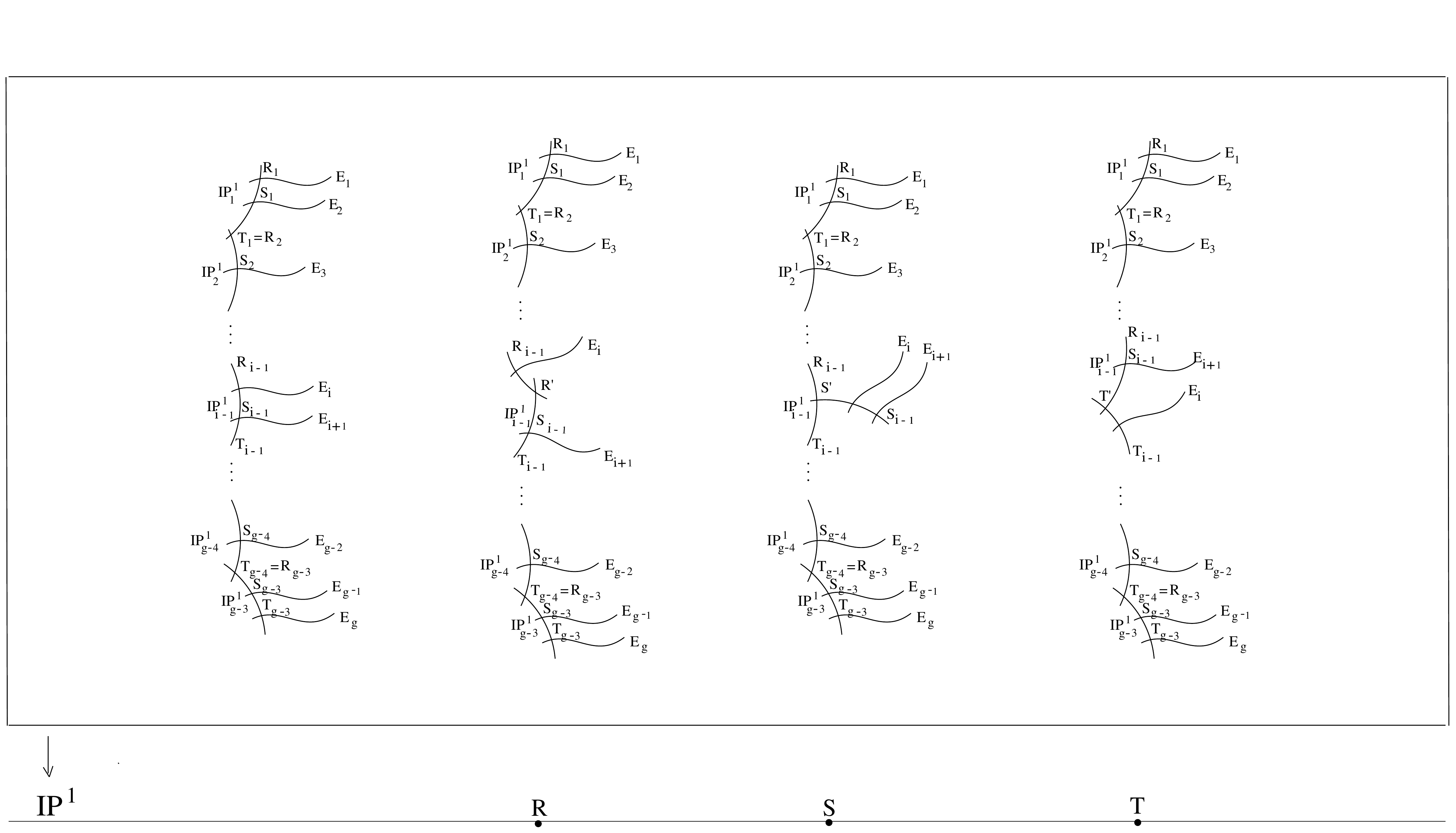}
	\caption{The family.}
	\label{fig:Familia}
\end{figure}
%C:/Users/Gabriel/Downloads/Familia.pdf
}
\end{subs}

\begin{prop}\label{ecuaciones}
Let $g\geq 5$ and let $D\subseteq \overline{M_{g}}$ be an effective divisor with class
\begin{center}
$D=a\lambda-a_0\delta_0-a_1\delta_1-\ldots-a_{[g/2]}\delta_{[g/2]}$.
\end{center}
If $[\pi_{i}]^{*}D=0$, for every $2\leq i\leq [g/2]$, then 
\begin{center}
$a_l=(l(g-l)/(g-1))a_1$, for every $2\leq l\leq [g/2]$.
\end{center}
\end{prop}
{\em Proof.} For every $i$, the degree of $(\delta_{0})_{\pi_i}$ is $0$, because each fiber of $\pi_i$ contains only disconnecting nodes. On the other hand, for every fiber $F$ of $\pi_i$, each section of $H^{0}(\omega_{\pi_i}\big|_{F})$ vanishes at each $\mathbb{P}^{1}$. 

Furthermore, we have that $H^{0}(\omega_{\pi_i}\big|_{E})=H^{0}(\omega_{E}(A))$ for every elliptic component $E$ of $F$, where $A$ is the node of $F$ lying on $E$.
The upshot is that
\begin{center}
$H^{0}(\omega_{\pi_i}\big|_{F})=\bigoplus_{E} H^{0}(\omega_{E})$,
\end{center}
for every fiber $F$ of $\pi_i$. Thus, $\pi_{i*}(\omega_{\pi_i})$ is trivial and hence deg$((\lambda)_{\pi_i})=0$.

Assume $i\geq 4$. By the construction of $\pi_i$, we have that
\begin{equation*}
\text{deg}((\delta_{1})_{\pi_i})=\widetilde{\Delta}^{2}+((\mathbb{P}^{1}_{i-1}\times\{S_{i-1}\})^{\widetilde{}})^{2}=\Delta^{2}-3+((\mathbb{P}^{1}_{i-1}\times\{S_{i-1}\})^{2}-1)=-2
\end{equation*}
%&=2-2(0)-3+(0-1)=-2,
where $\widetilde{\Delta}$ and $(\mathbb{P}^{1}_{i-1}\times\{S_{i-1}\})^{\widetilde{}}$ are the strict transforms of the diagonal $\Delta$ and $\mathbb{P}^{1}_{i-1}\times\{S_{i-1}\}$ in the blow up $(\mathbb{P}^{1}_{i-1}\times \mathbb{P}^{1}_{i-1})^{\widetilde{}}$ of $\mathbb{P}^{1}_{i-1}\times \mathbb{P}^{1}_{i-1}$ at the points $(R_{i-1},R_{i-1}),(S_{i-1},S_{i-1})$ and $(T_{i-1},T_{i-1})$.

On the other hand, we have deg$((\delta_{2})_{\pi_i})=1$, as the fiber of $\pi_i$ over $S$ has a disconnecting node $S'$ such that the closure of one of the connected components of $\pi_{i}^{-1}(S)-\{S'\}$ has genus $2$ and the total space of $\pi_i$ is smooth at $S'$ (the total space of $\pi_i$ is smooth at the point $S'$, as this point can be seen as a point of $(\mathbb{P}^{1}_{i-1}\times \mathbb{P}^{1}_{i-1})^{\widetilde{}}$, which is a smooth surface).

For $3\leq l\leq i-2$, we have deg$((\delta_{l})_{\pi_i})=0$, as the family is locally trivial around $\mathbb{P}^{1}\times \{Q\}\subseteq \cal{X}$ for every $Q$ which is a node of $X_i$ or $Y_i$.

Now, we will compute deg$((\delta_{i-1})_{\pi_i})$. Notice that for every fiber $F$ of $\pi_i$, the closure of one of the connected components of $F-\{R_{i-1}\}$ has genus $i-1$. If $g$ is even and $i=g/2$, then for each fiber $F$, the closure of one of the connected components of $F-\{T_{i-1}\}$ has genus $i-1$ and hence
\begin{equation*}
\begin{split}
\text{deg}((\delta_{i-1})_{\pi_i})&=((\mathbb{P}^{1}_{i-1}\times\{R_{i-1}\})^{\widetilde{}})^{2}+((\mathbb{P}^{1}_{i-1}\times\{T_{i-1}\})^{\widetilde{}})^{2}\\
&=((\mathbb{P}^{1}_{i-1}\times\{R_{i-1}\})^{2}-1)+((\mathbb{P}^{1}_{i-1}\times\{T_{i-1}\})^{2}-1)\\
&=-2.
\end{split}
\end{equation*} 
Otherwise, 
\begin{equation*}
\text{deg}((\delta_{i-1})_{\pi_i})=((\mathbb{P}^{1}_{i-1}\times\{R_{i-1}\})^{\widetilde{}})^{2}=-1.
\end{equation*} 
To compute deg$((\delta_{i})_{\pi_i})$, first notice that the fiber of $\pi_i$ over $R$ has a disconnecting node $R'$ such that the closure of one of the connected components of $\pi_{i}^{-1}(R)-\{R'\}$ has genus $i$ and the total space of $\pi_i$ is smooth at $R'$, and the same holds for the fiber of $\pi_i$ at $T$. Now, if $g$ is odd and $i=(g-1)/2$, then for each fiber $F$, the closure of one of the connected components of $F-\{T_{i-1}\}$ has genus $i$ and hence
\begin{equation*}
\text{deg}((\delta_{i})_{\pi_i})=2+((\mathbb{P}^{1}_{i-1}\times\{T_{i-1}\})^{\widetilde{}})^{2}=1.
\end{equation*} 
Otherwise, deg$((\delta_{i})_{\pi_i})=2$. Finally, if $i\leq [g/2]-1$, then
\begin{equation*}
\text{deg}((\delta_{i+1})_{\pi_i})=((\mathbb{P}^{1}_{i-1}\times\{T_{i-1}\})^{\widetilde{}})^{2}=-1
\end{equation*}
and deg$((\delta_{l})_{\pi_i})=0$, if $i+2\leq l\leq [g/2]$. Now, by hypothesis we have $[\pi_{i}]^{*}D=0$ for every $2\leq i\leq [g/2]$. So, using the degrees we have computed, we get
\begin{center}
$2a_1-a_2+a_{i-1}-2a_i+a_{i+1}=0$ for every $2\leq i< [g/2]$,
\end{center}
\begin{center}
$2a_1-a_2+2a_{i-1}-2a_i=0$, if $g$ is even and $i=g/2$, and
\end{center}
\begin{center}
$2a_1-a_2+a_{i-1}-a_i=0$, if $g$ is odd and $i=(g-1)/2$.
\end{center}
For $i=2,3$, analogously, we get the same equations. Now, solving the system of $[g/2]-1$ equations, we get that
\begin{center}
$a_l=(l(g-l)/(g-1))a_1$, for every $2\leq l\leq [g/2]$.
\end{center}
   \hfill $\Box$

\begin{thm}\label{teoremaprincipal}
Let $g\geq 5$ be an odd integer and let $\overline{S^{2}W}\subseteq \overline{M_g}$ be the effective divisor which is defined as
the closure of the locus of smooth curves $C$ with a pair of points $(P,Q)$ satisfying that $Q$ is a ramification point of the linear system $H^{0}(\omega_{C}(-P))$ with ramification weight at least $3$. 

Write the class of $\overline{S^{2}W}$ in $\mathrm{Pic}_{\mathrm{fun}}(\overline{M_{g}})\otimes \mathbb{Q}$ in the form
\begin{center}
$\overline{S^{2}W}=a\lambda-a_0\delta_0-a_1\delta_1-\ldots-a_{[g/2]}\delta_{[g/2]}$.
\end{center}
Then 
\begin{center}
$a=9g^5-51g^4+129g^3-207g^2+174g-54$, and 
\end{center}
\begin{center}
$a_i=(i(g-i)/(g-1))a_1$ for every $2\leq i\leq [g/2]$.
\end{center}
\end{thm}
{\em Proof.} We have already computed the coefficient $a$ (see Subsection \ref{coeffa}). On the other hand, it follows from Proposition \ref{propflag} that $[\pi_{i}]^{*}\overline{S^{2}W}=0$ for every $2\leq i\leq [g/2]$. Then, by Proposition \ref{ecuaciones}, $a_i=(i(g-i)/(g-1))a_1$ for every $2\leq i\leq [g/2]$.   \hfill $\Box$

\textbf{Acknowledgements.} The author would like to thank Eduardo Esteves for several helpful discussions. Also, we acknowledge the use of Singular [S] for some of the computations.

\bibliographystyle{alpha}

\end{document}